\documentclass[letterpaper, 10 pt, conference]{ieeeconf}
\usepackage{bbm}
\usepackage{amsmath,amssymb,amsfonts}
\usepackage{mathrsfs}
\usepackage{xcolor}
\usepackage{graphicx}
\usepackage{amsfonts}
\usepackage[english]{babel}
\usepackage{epstopdf}
\usepackage{textcomp}
\usepackage[linesnumbered,ruled,vlined]{algorithm2e}

\newtheorem{theorem}{Theorem}

\newtheorem{lemma}{Lemma}

\newtheorem{lemma*}{Lemma}

\renewcommand{\hat}{\widehat}

\newcommand{\cT}{\mathcal{T}}

\newcommand{\nn}{\nonumber}    

\newcommand{\defeq}{\buildrel\triangle\over =}
\newcommand{\sq}[1]{\left[#1\right]}
\newcommand{\cm}[1]{\left(#1\right)}
\newcommand{\bm}[1]{\left\{#1\right\}}

\newcommand{\figref}[1]{Fig.~\ref{#1}}

\newcommand{\mP}{\mathbb{P}}
\newcommand{\mE}{\mathbb{E}}
\newcommand{\mR}{\mathbb{R}}

\newcommand{\cX}{\mathcal{X}}

\newcommand{\cA}{\mathcal{A}}

\newcommand{\cK}{\mathcal{K}}

\newcommand{\cP}{\mathcal{P}}

\newcommand{\cH}{\mathcal{H}}

\newcommand{\tgamma}{\tilde{\gamma}}
\newcommand{\tsigma}{\tilde{\sigma}}
\newcommand{\ttheta}{\tilde{\theta}}

\newcommand{\pih}{\hat{\pi}}
\newcommand{\muh}{\hat{\mu}}

\newcommand{\Fh}{\hat{F}}

\newcommand{\sigf}{\sigma^f}

\newcommand{\siglf}{{\sigma^l,\sigma^f}}
\newcommand{\tgammatf}{{\tilde{\gamma}_t}^f}

\newcommand{\tgammatl}{{\tilde{\gamma}_t}^l}

\usepackage{xspace}
\usepackage[acronym,nogroupskip,nonumberlist,nopostdot]{glossaries}
\glsdisablehyper

\newacronym{mdp}{MDP}{Markov decision process}
\newacronym{ne}{NE}{Nash equillibrium}
\newacronym{mfe}{MFE}{mean-field equillibrium}
\newacronym{mpe}{MPE}{Markov perfect equillibrium}
\newacronym{mfg}{MFG}{mean-field game}
\newacronym{rl}{RL}{reinforcement learning}
\newacronym{marl}{MARL}{multi-agent reinforcement learning}
\newacronym{iot}{IoT}{Internet of Things}
\newacronym{ssg}{SSG}{Stackelberg security game}
\newacronym{pse}{PSE}{perfect Stackelberg equilibrium}
\newacronym{mpse}{MPSE}{Markov \gls{pse}}
\newacronym{pomdp}{POMDP}{partially observed \gls{mdp}}
\newacronym{mc}{MC}{Monte Carlo}
\newacronym{pbe}{PBE}{perfect Bayesian equilibrium}
\newacronym{gmfg}{GMFG}{graphon mean-field game}
\newacronym{gmfe}{GMFE}{graphon mean-field equillibrium}
\newacronym{spe}{SPE}{sub-game perfect equillibrium}
\newacronym{mkv}{MKV}{McKean-Vlasov}

\newcommand{\mdp}{\gls{mdp}\xspace}

\newcommand{\rl}{\gls{rl}\xspace}

\newcommand{\pse}{\gls{pse}\xspace}
\newcommand{\mpse}{\gls{mpse}\xspace}

\newcommand{\pbe}{\gls{pbe}\xspace}

\allowdisplaybreaks
\title{\LARGE \bf Model-free Reinforcement Learning for \\Stochastic Stackelberg Security Games}
\author{Rajesh K Mishra,
Deepanshu Vasal,
	and Sriram Vishwanath}
\begin{document}

	\maketitle

	\begin{abstract}
		In this paper, we consider a sequential stochastic Stackelberg game with two players, a leader and a follower. The follower has access to the state of the system while the leader does not. Assuming that the players act in their respective best interests, the follower's strategy is to play the best response to the leader's strategy. In such a scenario, the leader has the advantage of committing to a policy which maximizes its own returns given the knowledge that the follower is going to play the best response to its policy. Thus, both players converge to a pair of policies that form the Stackelberg equilibrium of the game. Recently,~\cite{Vasal} provided a sequential decomposition algorithm to compute the Stackelberg equilibrium for such games which allow for the computation of Markovian equilibrium policies in linear time as opposed to double exponential, as before. In this paper, we extend the idea to a \mdp whose dynamics are not known to the players, to propose a \rl algorithm based on Expected Sarsa that learns the Stackelberg equilibrium policy by simulating a model of the \mdp. We use \emph{particle filters} to estimate the belief update for a common agent which computes the optimal policy based on the information which is common to both the players. We present a security game example to illustrate the policy learned by our algorithm.
	
	\end{abstract}

	\section{Introduction}
Stackelberg games are a very useful tool to model strategic interactions where there is a dominant player called the leader who commits to a policy and a follower that \emph{observes} the leader's policy and plays its best response to it. With the knowledge that the follower will play the best response, the leader can devise an optimum strategy that maximizes its own rewards. Such games, also termed as Stackelberg security games, have become very popular in the recent decade in building and analyzing real world security systems in areas such as airports, seaports, and wildlife parks~\cite{fang2015security}. In such games, the leader has an inherent advantage that enables it to commit to a strategy that benefits it irrespective of the strategy followed by the follower given that it is rational and will play the best response. Such one shot equilibrium games have also proven useful for economic firms to analyze markets and take competitive risks.

Most of the prior work consider single shot Bayesian game models where the leader and the follower interact only once. However, most practical systems entail a periodic interaction between the leader and the follower. Computing Stackelberg equilibria for such stochastic games was unknown. Recently~\cite{Vasal} presented a sequential decomposition algorithm to compute \mpse of such games. Solving a dynamic stochastic Stackelberg game when the follower has a private Markovian state is computationally challenging. This is because, unlike other games, dynamic games of asymmetric information have players' strategies that are coupled across time. Since strategy of a player is a mapping from each history of the game which grows exponentially with time, the space of strategies of the players is double exponential, rendering such problems intractable. Recently, there has been results on sequential decomposition of certain classes of games of asymmetric information~(\cite{Vasal2016,Vasal2017,Vasal}). In repeated Stackelberg security games, there have been other approaches to mitigate this issue. For instance,~\cite{Kar2015} considers a repeated Stackelberg game and uses a new human behavior model to study such games. \cite{Balcan2015} consider a learning theoretic approach to study a repeated Stackelberg game between attacker and defender where they use regret analysis to learn attacker types, and show sub-linear regret for both complete and partial information models. There also have been efforts to develop \rl algorithms to study models with asymmetric information in~\cite{Kononen2003} for repeated Bayesian Stackelberg games.

Kalman filters have been widely used for Gaussian state space modelling but the advent of sequential monte carlo methods can be attributed to certain applications where non-Gaussian state space modelling was required~\cite{kitagawa1996monte}. Particle filters are sequential monte carlo filters that approximate the belief on a state, in other words, the probability of a system of being in a particular state, from an empirical distribution based on the observed history when model dynamics are unknown~\cite{Douc2005,Crisan2002}. It is widely popular in applications like robotics for localizations and fault dynamics, where most of the times the environment is non-Gaussian and needs to learned based on collected samples and observations~\cite{rekleitis2004particle}. These methods utilize a $K$ number of random samples or \emph{particles}, where $K$ is large, to represent the posterior probability of the state based on the observations.

In this paper, we propose an \rl algorithm with particle filters to learn the Stackelberg equilibrium strategies for a rational, leader and follower, when the players are unaware of the dynamics of the game. The algorithm learns the $Q$ values using the Expected Sarsa and then solves a fixed point equation for the follower and a maximization equation for the leader to converge upon the strategies. We use a common agent approach~\cite{Nayyar2013a} wherein a fictitious common agent has access to the common information that is shared between the players and uses it to compute a belief on the private state. It then solves for the optimal policies using the \rl algorithm and updates the belief using the particle filters. The use of \rl algorithm with particle filters for solving for Stackelberg equilibrium is novel. We illustrate our algorithm by determining the strategies for a realistic security game example to show that the algorithm derived, coincides with the optimal strategy that was obtained in~\cite{Vasal}, where it was assumed that the dynamics of the game were known.

The paper is structured as follows. The model is presented in Section~\ref{sec:Model} followed by a discussion on the sequential decomposition algorithm in Section~\ref{sec:Preliminaries}. We present our proposed algorithm in Section~\ref{sec:Algorithm} and prove convergence to the Stackelberg equilibrium in Section~\ref{sec:Conclusion}. In Section~\ref{sec:Example}, we provide an example to showcase our results and conclude in Section~\ref{sec:Conclusion}.

\subsection{Notation}

We use uppercase letters for random variables and lowercase for their realizations. For any variable, subscripts represent time indices and superscripts represent player identities, `$l$' for the leader and `$f$' for the follower. For any finite set $\mathcal{S}$, $\mathcal{P}(\mathcal{S})$ represents space of probability measures on $\mathcal{S}$ and $|\mathcal{S}|$ represents its cardinality. We denote by $P^{\sigma}$ (or $E^{\sigma}$) the probability measure generated by (or expectation with respect to) strategy profile $\sigma$ and the space for all such strategies as $\mathcal{K}^\sigma$. We denote the set of real numbers by $\mathbb{R}$. All equalities and inequalities involving random variables are to be interpreted in an \emph{a.s.} sense. 
	\section{Model}
	\label{sec:Model}
		
Consider a stochastic Stackelberg game over a time horizon $\cT$ with perfect recall between two players: a leader and a follower. The actions and the states are defined over finite sets $\cA$ and $\cX $ respectively. $\cA^l$ denotes the action set of the leader while $\cA^f$ represents that of the follower. The state $x_t$ evolves as the following controlled Markov process
\begin{align}
	\cP\left(x_t|a_{1:t-1},x_{1:t-1}\right) &=  \tau\left(x_t|a_{t-1},x_{t-1}\right), 
\end{align}
where $a_t = \cm{a_t^l,a_t^f}$ are the actions of the leader and the follower. The action history $a_{1:t-1}$ is shared between both the players as common information. However, only the follower has access to the private information $x_{1:t}$. The leader observes the common information $a_{1:t-1}$ and takes action $a_t^l\in \cA^l$ while the follower observes the private information $x_{1:t}$ and the common information $a_{1:t-1}$ and takes action $a_t^f\in \cA^f$.

	\section{Preliminaries}
	\label{sec:Preliminaries}
		
In this section, we discuss the sequential decomposition algorithm which used to compute the Stackelberg equilibrium policies for the case when the model is known to the players.

\subsection{Stackelberg Equilibrium}
Given a strategy profile $\sigma^l$ for the leader, the follower maximizes its own total discounted expected rewards over a finite horizon $T$ as
\begin{align}
	\max_\sigf\mE^\siglf \sq{\sum_{t=1}^T \delta^{t-1}R_t^f\left(X_t,A_t\right)}.
\end{align}
Now, if we denote $\Lambda^f(\sigma^l)$ to be the set of all optimizing strategies for the follower, given a strategy $\sigma^l$ of the leader, we get
\begin{align}
	\Lambda^f(\sigma^l) = \arg\max_{\sigma^f} \mE^{\sigma^l,\sigma^f} \left[ \sum_{t=1}^T\delta^{t-1} R_t^f\left(X_t,A_t\right) \right].
\end{align}
With the information that follower is going to play the best response to its strategy, the leader tries to maximize its own total expected discounted rewards by finding an optimal strategy. In other words, knowing the follower's counter strategy to be $\Lambda^f(\sigma^l)$, the leader tries to solve the following maximization equation to compute its own strategy:
\begin{align}
\label{eq:Stck_eq}
	\tsigma^l \in\max_{ \sigma^l} \mE^{\sigma^l,\Lambda^f(\sigma^l)} \left\{ \sum_{t=1}^T\delta^{t-1} R_t^l(X_t,A_t) \right\}.
\end{align}
Both leader and follower successively play this game and converge upon a pair of strategies $(\tsigma^l,\tsigma^f)$ constituting a Stackelberg equilibrium where $\tsigma^f \in \Lambda^f\left(\tsigma^l\right)$.

\subsection{Perfect Stackelberg equilibrium}
\label{sec:PBSE}

The notion of \pse, as introduced in~\cite{Vasal2019a}, is in line with the \pbe discussed in~\cite{Fudenberg1991}. In the context of the game, the leader and the follower play their actions $a_t^l\sim \tsigma^l_t\left(\cdot|a_{1:t-1}\right)$ and $a_t^f\sim \tsigma^f_t(\cdot|a_{1:t-1},x_t)$ with $\tsigma = \left(\tsigma^f,\tsigma^l\right)$ constituting the Stackelberg equilibrium.

Assume $(\tsigma,\mu)$ is a \pse of the game, where $\mu = (\mu_t)_{t\in[T]}$ and $\tsigma = \left(\tsigma^f,\tsigma^l\right)$ with $\tsigma^f \in \Lambda_t^f(\tsigma^l)\ \forall \ t\in[T]$. For any $t$ and action history $a_{1:t-1}$, $\mu_t[a_{1:t-1}]\in\cP(\cX)$ is the equilibrium belief on the current follower's state $x_t$, i.e. $\mu_t[a_{1:t-1}](x_t)=\mP^{\tsigma}(x_t|a_{1:t-1})$. For any given $\sigma^l$, $\Lambda_t^f(\sigma^l)$ be defined as~$\forall h_t^f$,
\begin{align}
 	\label{Eqn:BRf}
 	\Lambda_t^f(\sigma^l) \defeq \arg\max_{\sigma^f} \mE^{\sigma^l,\sigma^f,\mu_t} \left[ \sum_{n=t}^T \delta^{n-t}R_n^l(X_n,A_n) |h_t^f\right]
 \end{align}
and
\begin{align}
	\label{Eqn:sigmaL}
 	\tsigma^l \in\max_{ \sigma^l} \mE^{\sigma^l,\Lambda^f(\sigma^l),\mu_t} \left[ \sum_{n=t}^T \delta^{n-t}R_n^l(X_n,A_n) |h_t^c\right].
\end{align}
where $h_t^f =(a_{1:t-1},x_{1:t})$ and $h_t^c = a_{1:t-1}$.

For the case of \pse the sum of expected returns can be expressed as 
\begin{align}
	\label{Eqn:j_func1}
	&J^{f,\sigma^f,\sigma^l,\pi}_t=\nn\\
	&\mathbb{E}^{\sigma^f_{t:T}, \sigma_{t:T}^l,\pi_t}\left[\sum_{k=t}^{T}\delta^{k-t}R_t^f\left(X_k,A_k\right)\vert x_{1:t},a_{1:t-1}\right],
\end{align}
\begin{align}
	\label{Eqn:j_func2}
	J^{l,\sigma^f,\sigma^l,\pi}_t=\mathbb{E}^{\sigma^l_{t:T},\sigma_{t:T}^f,\pi_t}\left[\sum_{k=t}^{T}\delta^{k-t}R_t^l\left(X_k,A_k\right)\vert a_{1:t-1}\right].
\end{align}

\subsection{Common agent approach}

We use the common agent approach in line with the common information approach that was used in~\cite{Nayyar2013a}. An arbitrary common agent with access to the common information $a_{1:t-1}$ generates prescription functions $\gamma_t = (\gamma_t^l,\gamma_t^f)$. The prescription functions $\gamma_t^l\in \cP(\cA^i)$ and $\gamma_t^f:\cX^t\to  \cP(\cA^i)$ are used by the leader and the follower to generate their actions as $a_t^l \sim\gamma_t^l(\cdot)$ and $a_t^f \sim\gamma_t^f(\cdot|x_{1:t})$. 

We denote the prescription functions as $\gamma_t = (\gamma_t^l,\gamma_t^f)$ as $\theta_t[\pi_t]$ with $\pi_t(x_t) = P^{\theta}(x_t|a_{1:t-1})$. In other words, the Markovian common agent computes the prescription functions as function of the belief state $\pi_t$, which it derives from the common information $a_{1:t-1}$. The belief state $\pi_t$ is given as
\begin{align}
	\pi_t = \mP^{\ttheta}\cm{x_t\vert a_{1:t}}
\end{align}
which denotes the distribution of the state conditioned on the observed action history. The follower uses $\gamma_t^f$ to operates on its current private state $x_t$ to produce its action $a_t^f$, i.e. $\gamma_t^f : \cX\to \cP(\cA^f)$ and $a_t^f \sim\gamma_t^f(\cdot|x_t)$. while the leader uses $\gamma_t^l$ to produce its action $a_t^l$ as $\gamma_t^l \in \cP(\cA^f)$ and $a_t^l \sim\gamma_t^f(\cdot)$.

In order to track the belief state $\pi_t$ we derive a recursive equation given a policy $\ttheta_t$. 

\begin{lemma} For any given policy of type $\theta$, there exists update functions $F$, independent of $\theta$, such that
\label{Lemma:belief_state}
	\begin{align}
	\label{eq:F_update}
		\pi_{t+1} = F(\pi_t,\gamma_t^f,a_t) 
	\end{align}
\end{lemma}
which can be elaborated using the Bayes' theorem as
\begin{align}
\label{Eqn:pit_recursive}
	\pi_{t+1}\left(x_{t+1}\right) &= \mP^\theta\left(x_{t+1}\vert a_{1:t}\right)\nonumber\\
								  &=\frac{\sum_{x_t\in\cX}\pi_t\left(x_t\right)\gamma_t^f\left(a_t^f\vert x_t\right)Q_x\left(x_{t+1}\vert x_t, a_t\right)}{\sum_{x_t\in\cX}\pi_t\left(x_t\right)\gamma_t^f\left(a_t^f\vert x_t\right)}
\end{align}

In summary, at time $t$, the common agent observes the action history $a_{t-1}$ and generates the optimal policy $\ttheta$ as a function of the belief $\pi_t$. Corresponding to each belief $\pi_t\in\Pi$, the generated prescription function $\tgamma_t=\cm{\tgamma_t^f,\tgamma_t^l}$ specifies the actions to be taken by the players. The optimal policy and the actions thus generated are used to obtain the next belief $\pi_{t+1}$ using~\eqref{Eqn:pit_recursive}. 

\subsection{Particle Filters}
\label{sec:particle_filter}
The main challenge in computing the updated belief $\pi_{t+1}$ from current belief $\pi_t$ without knowledge of the transition function $\tau$.  Particle filters are sequential monte carlo filters that approximate the state distribution from an empirical distribution based on the observed history. These methods utilize a $K$ number of random particles to represent the posterior probability of the state based on the observations.

These filters approximates the belief state $\pi_t$ by a set of $K$ sampled points from the state space $x_t\in\cX$, updated in a sequential manner at every observation point $a_t$, which also serves as an action in our case, through a selection procedure to establish the truthfulness of the belief based on the observation. In other words, the belief using a particle filter could be expressed as,
  \begin{align}
    \pih_t=\frac{1}{K}\sum_{i=1}^Kf\cm{x_t}.
  \end{align}
  The generic particle filter called the bootstrap filter that samples the states from a previous distribution and the resamples based on observations. It is estimated using an empirical distribution given as 
  \begin{align}
    \pih^{i}\cm{x_t} \defeq \sum_{j=1}^K\delta_{x_t}w_t^j,
  \end{align}
  where $\delta_{x_t}$ is a dirac delta function made up of $K$ particles $x_t^{1:K}$. The algorithm recursively consists of two steps, a transition step to sample $K$ particles from the current distribution and to obtain the samples corresponding to the next states for each of the sample according to the transition function. This is followed by a selection step where it is resampled according to the weights $w_t$ generated based on the observations. The algorithm can be summarized as below~\cite{Coquelin,Doshi2009,Elvira2017}.

  \begin{enumerate}
    \item Initialize, $t = 0$, $x_0^i\sim\pi_t\cm{\cdot}$, set $t=1$,
    \item For $t = 1,2,3,\ldots$,\newline 
    \textbf{Importance sampling:}
    \begin{enumerate}
      \item For $ i =1 \ldots K$, sample from the model, $\hat{x}_i\sim\mP\cm{x_{t+1}\vert x_t,a_t}$
      \item For $ i = 1\ldots K$, compute the weights in proportion to the chances of the next state with the current observation $a_t$
      \begin{align}
        w_t^i=\frac{\mP^\theta\cm{a_t,\hat{x}_i}}{\sum_{j=1}^K\mP^\theta\cm{a_t,\hat{x}_j}}
      \end{align}
    \end{enumerate}
    \textbf{Selection/ resampling}
    \begin{enumerate}
      \item Resample from the list $\cm{\hat{x}_{t+1}}$ with replacement according to the weights to get $\cm{x_{t+1}}$. This is done by choosing from indices $\bm{1\ldots K}$ according to the multinomial distribution $\cm{w_1\ldots w_K}$.
      \item The new belief state estimate is then,
      \begin{align}
        \pih\cm{x_t} \defeq \sum_{i=1}^K\delta_{x_{t+1}}
      \end{align}
    \end{enumerate}
  \end{enumerate}
  
  The multinomial distribution used for resampling is one of the simplest methods that was introduced in~\cite{Vapnik1999}. This methods redistributes the samples based on their corresponding weights. Other versions of the resampling method include the stratified sampling method that reduces variance~\cite{kitagawa1996monte,Douc2005}.

  The common agent employs a $2$ parallel particle filters that estimate $\mP^\theta\cm{x_{t+1}\vert a_{1:T}}$ for each of the player given the policy $\theta$. The particle filter as a module takes in the current belief vector, the corresponding policy and the observation vector. It uses the model to sample the next steps and then computes the posterior distribution.

\subsection{Algorithm for \mpse computation}
\label{sec:Result}

In this section, we discuss the sequential decomposition algorithm stated in~\cite{Vasal2019a} but using the belief estimated using the particle filters to compute the equilibrium strategies by solving equations~\eqref{Eqn:BRf} and~\eqref{Eqn:sigmaL}.
 
\subsubsection{Backward Recursion}

In a general \rl setting, with finite state space $\cX$ and finite action space $\cA$, we can express the $Q$-value for the state action pair $\left(x_t, a_t\right)$ for all $x_t\in\cX$, $a_t\in\cA$ in terms of the reward $R_t\left(x_t, a_t\right)$ as
\begin{align}
\label{Eqn:Bellman_optimality}
	Q_t\left(x_t, a_t\right)=R_t\left(x_t, a_t\right)+\delta\mE^\tau\left[V_{t+1}\left(X_{t+1}\right)\right]. 
\end{align} 
Here $Q_x\left(\cdot\vert x_t, a_t\right)$ represents the dynamics of the \mdp.  $V_{t+1}\left(\cdot\right)$ gives the value function at the future state and $\delta$ is the discounting parameter. Our algorithm is based on Expected Sarsa, which is a model free \rl algorithm that does not require the knowledge of model dynamics and produces a low variance estimate of the $Q$-value function.

Let us define an equilibrium generating function $\ttheta=(\ttheta^i_t)_{i\in\{l,f\},t\in[T]}$, where $\ttheta_t : \mathcal{P}(\cX) \to \left\{\cX \to \mathcal{P}(\cA^f) \right\} \times \mathcal{P}(\cA^l) $. In addition, we define value functions $\cm{V_t^f,V_t^l}$ and action value functions $(Q_t^l,Q_t^f)_{t\in \{ 1,2, \ldots T+1\}}$ , where 
\begin{align}
	V_t^f&: \mathcal{P}(\cX) \times \cX \to \mR\nn\\
	V_t^l&: \mathcal{P}(\cX) \to \mR\nn\\
	Q_t^f&: \mathcal{P}(\cX) \times \cX \times \cA \times \Gamma^f \to \mR\nn\\
	Q_t^l&: \mathcal{P}(\cX) \times \cA \times \Gamma^f \to \mR\nn
\end{align}

Here we describe the backward recursive algorithm that is used to compute the strategies using the value functions for each discrete value of the estimated belief state $\pih_t$.

\begin{enumerate}
	\item Initialize $\forall \pih_{T+1}\in \mathcal{P}(\cX), x_{T+1}\in \cX$,
	\begin{align}
		\label{eq:VT+1}
		V^f_{T+1}\cm{\pih_{T+1},x_{T+1}} &\defeq 0 \\
		V^l_{T+1}\cm{\pih_{T+1}} &\defeq 0 
	\end{align} 		
	\item For $t = T,T-1, \ldots 1, \ \forall \pih_t\in \mathcal{P}(\cX)$, let $\ttheta_t[\pih_t]$ be generated as follows. 
	
	$\forall \pih_t\in\Pi$, $x_t\in\cX$, $a_t\in\cA$ and $\tgamma^f_t\in K^{\gamma^f}$, compute $Q_t^l$ and $Q_t^f$ as
	\begin{align}
		\label{Eqn:Qdef1}
		Q_t^f\cm{\pih_t,x_t,a_t,\tgamma^f_t}&= \mE^\tau\Big[R_t^f\left(x_t,a_t\right) +\nn\\
		\delta V_{t+1}^f\Big(&\Fh^f\left(\pih_t, \tgamma^f_t, a_t\right), X_{t+1}\Big) \lvert x_t\Big]
	\end{align} 
	\begin{align}
		\label{Eqn:Qdef2}
		Q_t^l\cm{\pih_t,a_t,\tgamma^f_t}&= \mE^{\tau,\pih_t}\Big[R_t^l\left(X_t,a_t\right) +\nn\\
		\delta V_{t+1}^f&\left(\Fh^l\left(\pih_t, \tgamma^f_t, a_t\right)\right)\Big]
	\end{align} 
	where $\Fh\cm{\cdot}$ represents the particle filter function.

	Set $\tgamma_t = \ttheta_t[\pih_t]$, where $\tgamma_t =\cm{\tgamma_t^l,\tgamma_t^f}$ is the solution of the following fixed-point equation. 
	
	$\forall \gamma_t^l\in K^{\gamma^l}$, define $\Lambda(\gamma_t^l)$ as follows, $\forall x_t\in \cX$,
	\begin{align}
		\label{Eqn:fixed_point}
		&\Lambda_t^f(\gamma_t^l) =\Big\{\tgamma_t^f:\tgamma_t^f\in\nn\\
		&\arg\max_{\gamma^f_t\cm{\cdot\vert x_t}} \mE^{\gamma_t^f\cm{\cdot\vert x_t},\gamma_t^l}\sq{Q_t^f\cm{\pih_t,x_t,A_t,\tgamma^f_t}} \Big\}.
	\end{align}

	The solution to the strategy for the leader is compute using
	\begin{align}
	\label{Eqn:fixed_point2}
	\tgamma_t^l &\in \arg\max_{\gamma_t^l}\mE^{ \Lambda_t^f\left(\gamma_t^l\right){\gamma}^{l}_t} \sq{Q_t^l\cm{\pih_t,A_t,\tgamma^f_t}}
	\end{align}

	The solution to the above two steps generate $(\tgamma_t^l,\tgamma_t^f)$ as pair of equilibrium strategies.
	
	Now, we calculate the utility function values $\forall x_t\in\cX$,
	\begin{align}
	\label{Eqn:utility1}
		V^f_{t}\left(\pi_t,x_t\right) \defeq  &\;\mE^{\tilde{\gamma}^{f}_t\left(\cdot|x_t\right) \tilde{\gamma}^{l}_t,\, \pi_t}\left[{R}_t^f\left(x_t,A_t\right) + \right.\nn\\
								   			  &\left.\delta V_{t+1}^f \left(\Fh^f\left(\pi_t, \tilde{\gamma}^f_t, A_t\right), X_{t+1}\right)\big\lvert x_t \right]
	\end{align}
	and
	\begin{align}
	\label{Eqn:Utility2}
		V^l_{t}\left(\pi_t\right) \defeq  &\;\mE^{\tilde{\gamma}^{f}_t \tilde{\gamma}^{l}_t,\, \pi_t}\Big[ {R}_t^l\left(X_t,A_t\right) +\nn\\
							   &\delta V_{t+1}^l \left(\Fh^l\left(\pi_t, \tilde{\gamma}^f_t, A_t\right)\right) \Big].	
	\end{align}
\end{enumerate}

\subsubsection{Forward Recursion}

Based on $\ttheta$ defined in the backward recursion above, we now construct a set of strategies $\tsigma$ (through beliefs $\mu$) in a forward recursive way as follows. 
\begin{enumerate} 
	\item Initialize at time $t=1$, 
	\begin{align}
	\mu_1[\phi]\left(x_1\right) &:=  \tau\left(x_1\right). \label{eq:mu*def0}
	\end{align}
	
	\item For $t =1,2 \ldots T, \forall i =1,2, a_{1:t}\in \cH_{t+1}^c, x_{1:t} \in\cX^t$
	\begin{align}
	\label{Eqn:beta*def}
		\tsigma_{t}^{l}\left(a_{t}^l|a_{1:t-1}\right) &:= \ttheta_{t}^l\left[\mu_{t}\left[h_t^c\right]\right]\left(a^l_{t}\right) \\
		\tsigma_{t}^{f}\left(a_{t}^f|a_{1:t-1},x_{1:t}\right) &:= \ttheta_{t}^f\left[\mu_{t}\left[h_t^c\right]\right]\left(a^f_{t}|x_{t}\right)
	\end{align}
	\begin{align}
	\label{Eqn:mu*def}
		\mu_{t+1}\left[h_{t+1}^c\right] :=F\left(\mu_t\left[h_t^c\right], \ttheta_t^f\left[\mu_t\left[h_t^c\right]\right], a_t\right)
	\end{align}
\end{enumerate}
where $F$ is defined in \eqref{eq:F_update}.

	\section{Reinforcement Algorithm}
	\label{sec:Algorithm}
		In this section, we propose a model-free \rl algorithm based on Expected Sarsa~\cite{van2009theoretical} to compute the optimal \mpse strategy both for the leader and the follower. We follow the steps summarized in the sequential decomposition algorithm to derive these strategies. The algorithm iterates between the policy evaluation and policy improvement steps sequentially. At each instant $t$, and at each discretized belief state $\pih_t\in\Pi$, we compute $Q_t^i$ as in~\eqref{Eqn:Qdef1} and~\eqref{Eqn:Qdef2} using Expected Sarsa as a function of the equilibrium policy of the follower $\ttheta\sq{\pih_t}=\tgammatf$. This is essential in solving the fixed point equations for model free algorithms. Thereafter, we compute the optimal policies for the follower and the leader by solving the fixed point equation in~\eqref{Eqn:fixed_point} and the maximization in~\eqref{Eqn:fixed_point2} respectively. The value functions $V_t^f$ and $v_t^l$ are updated after the optimal policy $\ttheta_t$ is computed using the equations in~\eqref{Eqn:utility1} and~\eqref{Eqn:Utility2} to be used in the computation of $Q_t^i$ in the next iteration. 

\subsection{Policy Evaluation}

The $Q$-value functions $Q_{t}^f$ and $Q_{t}^l$ are computed from a sampled trajectory from the model bootstrapping it with the future value function $V_{t+1}^i$. It involves sampling from a model simulation of a generic player $i\in\left\{l,f\right\}$ which takes a current state $x_t\in\cX$ and actions $a_t\in\cA$, for every belief state $\pih_t$ and provides the next state $x_{t+1}$ and the corresponding reward $R_t^i(x_t, a_t)$. In the beginning of each iteration, the belief state estimate is updated to $\pih_{t+1}$ from the current estimated belief state $\pih_t$ at the equilibrium policy $\gamma_t^f$ by feeding the observation vector $a_t$ to the particle filter. It is computed separately for the follower and the leader as the model is not known to either of them. We use Expected Sarsa algorithm to update the $Q$-value for both the follower and the leader synchronously using the same simulated trajectory. Expected Sarsa follows the TD learning method where the update equation is given as
\begin{align}
	Q^f_t\cm{\pih_t,x_t,a_t,\tgamma_t^f} = \cm{1-\alpha}Q^f_t\cm{\pih_t,x_t,a_t,\tgamma_t^f}\nn\\ +\alpha \cm{R_t^f\cm{x_t,a_t}+\delta V_{t+1}^f\cm{\pih^f_{t+1},x_{t+1}}}
\end{align}
\begin{align}
	Q^l_t\cm{\pih_t,a_t,\tgamma_t^f} = \cm{1-\alpha}Q^l_t\cm{\pih_t,a_t,\tgamma_t^f}\nn\\ +\alpha \cm{\mE^{\pih_t}[R_t^l\cm{X_t,a_t}]+\delta V_{t+1}^l\cm{\pih^f_{t+1}}}
\end{align}
where $\left(x_t, a_t, x_{t+1}\right)$ are sampled from a simulated model and $\alpha$ is the learning parameter. We use linear interpolation to obtain $V_{t+1}^f\cm{\pih_{t+1},x_{t+1}}$, $V_{t+1}^l\cm{\pih_{t+1}}$ at a future state. The updated belief state is generated using the particle filter algorithm that was put forth in~\ref{sec:particle_filter}. It takes in the current belief state $\pih_{t}$, the follower's policy $\tgammatf$ and the observation $a_t$ and provides the empirical distribution $\pih_{t+1}$. 

\subsection{Policy Iteration}
\subsubsection{Follower Strategy}

The fixed point equation in \eqref{Eqn:fixed_point} is solved at the follower, for each discrete belief state $\pih_t\in\Pi$, to compute the follower's equilibrium strategy corresponding to each of the leader's strategy $\gamma_t^l\in\cK^{\gamma^l}$. Given the follower's $Q_t^f$-value function, the objective function for $x_t\in\cX$ with any strategy $\gamma_t^f$ can be expressed as

\begin{align}
	\label{Eqn:objective1}
	\tgammatf = \arg\max_{\gamma_t^f} \mE\sq{Q_{t}^f\cm{\pih_t, x_t, A_t,\gamma_t^f}}.
\end{align}

The expectation is taken over $A_t$ through the measure $\gamma_t^f\left(\cdot\vert x_t\right)\gamma_t^l\left(\cdot\right)$ representing the leader's strategy. It is solved by a continuous policy update in the direction of ascent of the gradient. This is achieved using any of the policy gradient approaches. In this paper, we use neural network based policy gradient method to compute the policies that maximize the gradient of the $Q_t^f$ value functions at any time $t$. Given that the the function $Q_t^f$ depends on the optimal policy, this process is repeated over iterations arriving at the required prescription function $\tgamma^f_t$. This is repeated at all the belief states $pih_t\in\Pi$ so that we get the final equilibrium function $\theta_t\left[\pih_t\right]$.

\subsubsection{Leader Strategy}

The leader, with the knowledge of the strategy set of the follower $\Lambda_t^f\left(\gamma_t^l\right)$, computed as a best response to its strategies, optimizes its own strategy by solving the maximization equation in \eqref{Eqn:fixed_point2} given as 
\begin{align}
	\label{Eqn:objective2}
	\tgammatl = \arg\max_{\gamma_t^l} \mE\sq{Q_{t}^f\cm{\pih_t, A_t, \gamma_t^f}}.
\end{align}
It is a known fact that this results in pure strategies. This implies that the solution could be easily computed as the greedy policy $\tgamma_t^l$ that has the highest $Q$ value for each belief $\pih_t$.

Finally, the equilibrium strategy is given by $\tilde{\gamma}_t = \theta_t[\pi_t]$, where $(\tgamma_t^l,\tgamma_t^f)$ and $\tgamma_t^f = BR_t^f\left(\gamma_t^l)\right)$. 

\begin{algorithm}[!tb]
	\label{alg:Optimal_Policy}
	\DontPrintSemicolon
	\SetAlgoLined
	\KwIn{$Q^f_{t}$, $Q^l_{t}\quad \forall \gamma^f_{t}\in{\cK}^{\gamma^f_{t}}$, $\pih_{t} \in \Pi$\newline
		$\theta_{t}\left[\pih_{t}\right]\quad \forall \pih_t\in\Pi$}
	\vspace{.1cm}
	\KwOut{$\ttheta$}
	\vspace{.1cm}
	Initialize: $V^f_{T+1}$, $V^l_{T+1}$\;
	\For{$t=T\ldots1$}
	{
		Evaluate: $Q_t^f = FollowerQ\cm{V_{t+1}^f}$\;
		Evaluate: $Q_t^l = LeaderQ\cm{V_{t+1}^l}$\;
		\For{$\pih_t\in\Pi$}
		{
			\For{$\gamma^l_t\in\Pi$}
			{
				Compute: $\Lambda\cm{\gamma_t^l} = PolicyF\cm{Q_t^f,\gamma_t^l}$\;
			}
			Compute: $\tgammatl = PolicyL\cm{Q_t^l}$\;
			$\tgammatf = \Lambda\cm{\tgammatl}$\;
			\For{$x_t\in\cX$}
			{
				$V_t^f\cm{\pih_t,x_t} = \mE^{\tgammatf,\tgammatl}\sq{Q_t^f\cm{\pih_t,x_t,A_t}}$\;				
			}
			$V_t^l = \mE^{\tgammatf,\tgammatl}\sq{Q_t^l\cm{\pih_t,A_t}}$\;
			$\ttheta\sq{\pih_t}=\tgammatf$, $\tgammatl$\;
		}
	}
	\KwResult{$\ttheta$}
	\caption{Optimal policy}
\end{algorithm}

\begin{algorithm}[!tb]
	\label{alg:Policy Evaluation}
	\DontPrintSemicolon
	\SetAlgoLined
	\KwIn{$V_{t+1}^f$, $V_{t+1}^l$, $\gamma_t^l$}
	\vspace{.1cm}
	\KwOut{$Q_t^f$, $Q_t^l$}
	\vspace{.1cm}
	\For{$\pi_t,\gamma_t^f\in\Pi\times{\cK}^{\gamma^f}$}
	{ 
		Initialize: $Q_t^f$, $Q_t^l$\;
		\For{$l=1\ldots L$}
		{
			\For{$a^f_t,a^l_t\in\cA^f\times\cA^l$}
			{
				$a_t=\cm{a_t^f,a_t^l}$\;
				\For{$x_t\in\cX$}
				{
					Sample: $x_{t+1}\sim \tau\cm{\cdot\vert x_t,a_t}$\;
					$\pih^f_{t+1} = PF\cm{\pi_t, \gamma_t^f, a_t}$\; 
					$T^f =  R_t^f\cm{x_t,a_t}+\delta V_{t+1}^f\cm{\pih^f_{t+1},x_{t+1}}$
					$R^l{x_t,a_t} = R_t^l\cm{x_t,a_t}$\;
					$Q^f_t\cm{\pi_t,x_t,a_t,\gamma_t^f} = \cm{1-\alpha}Q^f_t\cm{\pi_t,x_t,a_t,\gamma_t^f} + \alpha T^f$\;
				}
				$\pih^l_{t+1} = PF\cm{\pi_t, \gamma_t^f, a_t}$\;
				$T^l = \mE^{\pi_t}\sq{R^l\cm{X_t,a_t}} + V_{t+1}^l\cm{\pih_{t+1}}$\;
				$Q^l_t\cm{\pi_t, a_t,\gamma^f_t} = \cm{1-\alpha}Q^l_t\cm{\pi_t, a_t,\gamma^f_t} + \alpha T^l$\; 
			}
				
		} 
	
	}	
	\KwResult{$Q^f_t$, $Q^l_t$}
	\caption{Policy Evaluation ($Q$-value Estimation)}
\end{algorithm}

	\section{Convergence}
	\label{sec:Convergence}
		
In this section, we prove the convergence of the proposed \rl algorithm to the equilibrium strategy of the statistical Stackelberg equilibrium. We put forth the following theorem.
\begin{theorem}
	The policies $\tsigma=\cm{\tsigma^f,\tsigma^l}$ thus generated from the \rl algorithm using the particle filter does form a $\epsilon$-\mpse of the Stackelberg game i.e.
	For the follower, we show,
	\begin{align}
		\label{Eqn:follower_1}
		J^{f,\tsigma^f,\tsigma^l}_t + \epsilon_t \geq J^{f,\sigma^f,\tsigma^l}_t
	\end{align}
	and for the leader we show,
	\begin{align}
		\label{Eqn:leader_1}
		J^{l,\tsigma^f,\tsigma^l}_t + \epsilon_t \geq J^{f,\Lambda\cm{\sigma^l},\sigma^l}_t.
	\end{align}

	\begin{proof}
		The proof to~\eqref{Eqn:follower_1} has been shown in Theorem~\ref{Thm:main_theorem} while the proof to~\eqref{Eqn:leader_1} has been shown in Theorem~\ref{Thm:main_theorem1}.
	\end{proof}

\end{theorem}

\begin{lemma} \label{Lemma:V_pf}
	The worst case error in the value function if we use a particle filter to update the belief state can be expressed as
	\begin{align}
		\label{Eqn:lemma2_to_prove}
		\Vert V^i_t\left(\pih_t,x_t^i\right) - {V^\star}^i_t\left(\pi_t,x_t^i\right)\Vert \leq \epsilon_1
	\end{align}
	where $\pi_t = \mu_t[1:a_{1-1}]$ and $\pih_t = \muh_t[1:a_{1-1}]$ are the belief trajectories with and without the particle filters.

	\begin{proof}
		This could be proved by using the Chernoff-Hoeffding bounds for deviation between sum of independent random variables from their true expectation. We leave the complete proof for the online version~\cite{schmidt1995chernoff}.

		It is also worth noting that the the value $epsilon$ is a function of the length of the time horizon and the discount factor. Moreover, it goes down with the increase in the number of particles used by the particle filter to estimate the belief. In our paper, we invariably choose a large $K$.
	\end{proof}
\end{lemma}

\subsection{Convergence of the Follower}
We put forth two lemmas that is used to prove the theorems. The lemmas have not been proven in this section due to lack of space but could be referred from the online version. 
\begin{lemma} 
	\label{Lemma:V_optimal_returns1}
	$\forall t\in\sq{T}$, $\forall \muh_t$ and $x_t\in\cX$,
	\begin{align}
		V_t^f\cm{\pih_t, x_t}= J^{f,\tsigma^f,\tsigma^f,\pih}_t
	\end{align}
	where $\tsigma_t=\cm{\tsigma^f_t,\tsigma_t^l}$ is the optimal policy at time $t$, $\pih_t = \muh_t[1:a_{1-1}]$ and $J_{t}$ is the accumulated optimal returns from $t$ till $T$.  
\end{lemma}

\begin{lemma} \label{Lemma:V_is_better1}
	$\forall t\in\cT$,$\forall \cm{a_{1:t-1},x_{1:t}}\in\cH_t^f$, $\sigma_t^f$,
	\begin{align}
		\label{Eqn:lemma2_to_prove1}
		{V_t^\star}^f\left(\pi_t,x_t\right) \geq \mE^{\sigma_t^f,\tsigma_{t}^l,\pi_t}\left[{Q^\star}^f_{t}\left(\pi_t, X_t, A_t,\tgamma^f_{t}\right)\right]
	\end{align}
	where $\pi_t = \mu_t[1:a_{1-1}]$ and $V^\star$ and $Q^\star$ represent the value functions in the case where model is known.
\end{lemma}

\begin{theorem}\label{Thm:main_theorem}
	A strategy $(\tilde{\sigma})$ constructed from the above algorithm is an $\epsilon$-\mpse of the game. i.e
	\begin{align}
		\label{Eqn:follower_main}
		J^{f,\tsigma^f,\tsigma^l,\pih_t}_t + \epsilon_2 \geq J^{f,\sigma^f,\tsigma^l,\pi_t}_t
	\end{align}

	We prove it through the technique of mathematical induction and will use the results that were proved before in Lemma~\ref{Lemma:V_optimal_returns2} and Lemma~\ref{Lemma:V_is_better2}.

	\begin{proof}
	For the base case, we consider $t = T$. The expected sum of returns, when the player $i$ follows the equilibrium policy $\tsigma$ is given as
	\begin{align}
			J_T^{f,\tsigma^f,\tsigma^l,\pi} = {V^\star}_T^f\left(\pi_T,x_T\right),\\
			J_T^{f,\tsigma^f,\tsigma^l,\pih} = V_T^f\left(\pih_T,x_T\right),
		\end{align}
		which is true from Lemma~\ref{Lemma:V_optimal_returns1}.
	From Lemma~\ref{Lemma:V_pf}, 
	\begin{align}
		\label{Eqn:follower_thm1}
		{V^\star}_T^f\left(\pi_T,x_T\right) \leq  V_T^f\left(\pih_T,x_T\right) + \epsilon
	\end{align}
	for some chosen small $\epsilon_T$. 
	From Lemma~\ref{Lemma:V_is_better2} and~\eqref{Eqn:follower_thm1} we get 
	\begin{align}
		J^{f,\tsigma^f,\tsigma^l,\pih_T}_T + \epsilon_T \geq J^{f,\sigma^f,\tsigma^l,\pi_T}_T
	\end{align}
		Assuming that the condition in~\eqref{Eqn:follower_main} holds at $t = t+1$, we get,
		\begin{align}
			\label{Eqn:assumption1}
			J^{f,\tsigma^f,\tsigma^l,\pih_{t+1}}_{t+1} + \epsilon_{t+1} \geq J^{f,\sigma^f,\tsigma^l,\pi_t}_{t+1}
		\end{align}
		We need to prove that the expression in~\eqref{Eqn:follower_main} holds for $t=t$ as well. 

		\begin{align}
			&J^{f,\tsigma^f,\tsigma^l,\pih_{t}}_{t}\label{Eqn:thm_le1}\\
			&=V^f_t\cm{\pih_t,x_t}\label{Eqn:thm_follower2}\\
			&\geq{V^\star}^f_t\cm{\pi_t,x_t} - \epsilon_{t+1}\label{Eqn:thm_follower3}\\
			&\geq\mE^{\sigma_t^f,\tsigma_{t}^l,\pi_t}\left[{Q_t^\star}^f\left(\pi_t, x_t, A_t, \tgamma^f_{t}\right)\right]- \epsilon_{t+1}\label{Eqn:thm_follower4}\\
			&=\mE^{\sigma_t^f,\tsigma_{t}^l,\pi_t}\left[R^f(x_t,a_t) + \delta {V^\star}^f_{t+1}\cm{\pi_{t+1},x_{t+1}}\right]- \epsilon_{t+1}\label{Eqn:thm_follower5}\\
			&\geq\mE^{\sigma_t^f,\tsigma_{t}^l,\pi_t}\left[R^f(x_t,a_t) + \delta {V}^f_{t+1}\cm{\pih_{t+1},x_{t+1}}\right]\label{Eqn:thm_follower6}- \epsilon_{t+1}\\
			&=\mE^{\sigma_t^f,\tsigma_{t}^l,\pi_t}\left[R^f(x_t,a_t) + \delta J^{f,\tsigma^f,\tsigma^l,\pih_{t+1}}_{t+1}\right]- \epsilon_{t+1}\label{Eqn:thm_follower7}\\
			&\geq\mE^{\sigma_t^f,\tsigma_{t}^l,\pi_t}[R^f(x_t,a_t) + \delta J^{f,\sigma_t^f,\tsigma_{t}^l,\pi_{t+1}}_{t+1}- \cm{1+\delta}\epsilon_{t+1}]\label{Eqn:thm_follower8}\\
			&=J^{f,\sigma_t^f,\tsigma_{t}^l,\pi_{t}}_{t} - \cm{1+\delta}\epsilon_{t+1}\label{Eqn:thm_follower9}
		\end{align}
		Thus, $J^{f,\tsigma^f,\tsigma^l,\pih_{t}}_{t} + \epsilon_{t}\geq J^{f,\sigma_t^f,\tsigma_{t}^l,\pi_{t}}_{t}$.

		~\eqref{Eqn:thm_follower2} is from Lemma~\ref{Lemma:V_optimal_returns2} and~\eqref{Eqn:thm_follower3} is from Lemma~\ref{Lemma:V_pf}.~\eqref{Eqn:thm_follower4} and~\eqref{Eqn:thm_follower5} are from standard definitions while ~\eqref{Eqn:thm_follower6} is true because ${V^\star}_{t+1}^l$ is optimal.~\eqref{Eqn:thm_follower7} is from Lemma~\ref{Lemma:V_optimal_returns2} followed by using assumption made at $t=t+1$ in~\eqref{Eqn:thm_follower8}.
	\end{proof}
\end{theorem}
\subsection{Convergence of the Leader}
We put forth two lemmas that is used to prove the theorems. The lemmas have not been proven in this section due to lack of space but could be referred from the online version.
\begin{lemma} 
	\label{Lemma:V_optimal_returns2}
	$\forall t\in\sq{T}$, $\forall \muh_t$,
	\begin{align}
		V_t^l\cm{\pih_t}= J^{l,\Lambda\cm{\tsigma^l},\tsigma^l,\pih}_t
	\end{align}
	where $\tsigma_t$ is the optimal policy at time $t$, $\pih_t = \muh_t[1:a_{1-1}]$ and $J_{t}$ is the accumulated optimal returns from $t$ till $T$.  
\end{lemma}

\begin{lemma} \label{Lemma:V_is_better2}
	$\forall t\in\cT$,$\forall \cm{a_{1:t-1}}\in\cH_t^f$,$\sigma_t^l$,
	\begin{align}
		\label{Eqn:lemma2_to_prove2}
		{V_t^\star}^l\left(\pi_t\right) \geq \mE^{\Lambda\cm{\sigma_t^l},\sigma_{t}^l,\pi_t}\left[{Q^\star}^l_{t}\left(\pi_t, A_t,\tgamma^f_{t}\right)\right]
	\end{align}
	where $\pi_t = \mu_t[1:a_{1-1}]$ and $V^\star$ and $Q^\star$ represent the value functions in the case where model is known.
\end{lemma}

\begin{theorem}\label{Thm:main_theorem1}
	A strategy $(\tilde{\sigma})$ constructed from the above algorithm is an $\epsilon$-\mpse of the game. i.e
	\begin{align}
		\label{Eqn:leader_main}
		J^{l,\tsigma^f,\tsigma^l,\pih_t}_t + \epsilon_t \geq J^{l,\Lambda\cm{\sigma^l},\sigma^l,\pi_t}_t
	\end{align}

	We prove it through the technique of mathematical induction and will use the results that were proved before in Lemma~\ref{Lemma:V_optimal_returns2} and Lemma~\ref{Lemma:V_is_better2}.

	\begin{proof}
	For the base case, we consider $t = T$. The expected sum of returns, when the player $i$ follows the equilibrium policy $\tsigma$ is given as
	\begin{align}
			J_T^{l,\Lambda\cm{\tsigma^l},\tsigma^l,\pi} = {V^\star}_T^l\left(\pi_T\right),\\
			J_T^{l,\Lambda\cm{\tsigma^l},\tsigma^l,\pih} = V_T^l\left(\pih_T\right),
		\end{align}
		which is true from Lemma~\ref{Lemma:V_optimal_returns2}.
	From Lemma~\ref{Lemma:V_pf}, 
	\begin{align}
		\label{Eqn:leader_thm1}
		{V^\star}_T^l\left(\pi_T\right) \leq  V_T^l\left(\pih_T\right) + \epsilon_T
	\end{align}
	for some chosen small $\epsilon_T$.

	From Lemma~\ref{Lemma:V_is_better2} and~\eqref{Eqn:leader_thm1} we get 
	\begin{align}
		J^{l,\tsigma^f,\tsigma^l,\pih_T}_T + \epsilon_T \geq J^{l,\Lambda\cm{\sigma^f},\tsigma^l,\pi_T}_T
	\end{align}
		Assuming that the condition in~\eqref{Eqn:leader_main} holds at $t = t+1$, we get,
		\begin{align}
			\label{Eqn:assumption2}
			J^{l,\tsigma^f,\tsigma^l,\pih_{t+1}}_{t+1} + \epsilon_{t+1} \geq J^{l,\Lambda\cm{\sigma^f},\tsigma^l,\pi_t}_{t+1}
		\end{align}
		We need to prove that the expression in~\eqref{Eqn:leader_main} holds for $t=t$ as well. 

		\begin{align}
			&J^{l,\tsigma^f,\tsigma^l,\pih_{t}}_{t}\label{Eqn:thm_leader1}\\
			&=V^l_t\cm{\pih_t}\label{Eqn:thm_leader2}\\
			&\geq{V^\star}^l_t\cm{\pi_t} - \epsilon_{t+1}\label{Eqn:thm_leader3}\\
			&\geq\mE^{\Lambda\cm{\sigma_t^l},\sigma_{t}^l,\pi_t}\left[{Q_t^\star}^l\left(\pi_t, x_t, A_t, \tgamma^f_{t}\right)\right]- \epsilon_{t+1}\label{Eqn:thm_leader4}\\
			&=\mE^{\Lambda\cm{\sigma_t^l},\sigma_{t}^l,\pi_t}\left[R^l(x_t,a_t) + \delta {V^\star}^l_{t+1}\cm{\pi_t}\right]- \epsilon_{t+1}\label{Eqn:thm_leader5}\\
			&\geq\mE^{\Lambda\cm{\sigma_t^l},\sigma_{t}^l,\pi_t}\left[R^l(x_t,a_t) + \delta {V}^l_{t+1}\cm{\pih_t}\right]\label{Eqn:thm_leader6}- \epsilon_{t+1}\\
			&=\mE^{\Lambda\cm{\sigma_t^l},\sigma_{t}^l,\pi_t}\left[R^l(x_t,a_t) + \delta J^{l,\tsigma^f,\tsigma^l,\pih_{t+1}}_{t+1}\right]- \epsilon_{t+1}\label{Eqn:thm_leader7}\\
			&\geq\mE^{\Lambda\cm{\sigma_t^l},\sigma_{t}^l,\pi_t}[R^l(x_t,a_t) + \delta J^{l,\Lambda\cm{\sigma^l},\sigma^l,\pi_{t+1}}_{t+1} \nn\\
			&- \cm{1+\delta}\epsilon_{t+1}]\label{Eqn:thm_leader8}\\
			&=J^{l,\Lambda\cm{\sigma^l},\sigma^l,\pi_{t}}_{t} - \cm{1+\delta}\epsilon_{t+1}\label{Eqn:thm_leader9}
		\end{align}
		Thus, $J^{l,\tsigma^f,\tsigma^l,\pih_{t}}_{t} + \epsilon_{t}\geq J^{l,\Lambda\cm{\sigma^l},\sigma^l,\pi_{t}}_{t}$.

		~\eqref{Eqn:thm_leader2} is from Lemma~\ref{Lemma:V_optimal_returns2} and~\eqref{Eqn:thm_leader3} is from Lemma~\ref{Lemma:V_pf}.~\eqref{Eqn:thm_leader4} and~\eqref{Eqn:thm_leader5} are from standard definitions while ~\eqref{Eqn:thm_leader6} is true because ${V^\star}_{t+1}^l$ is optimal.~\eqref{Eqn:thm_leader7} is from Lemma~\ref{Lemma:V_optimal_returns2} followed by using assumption made at $t=t+1$ in~\eqref{Eqn:thm_leader8}.
	\end{proof}

\end{theorem}

	\section{Numerical Example}
	\label{sec:Example}
		
In this section, we consider an example of a repeated Stackelberg security game~\cite{Vasal} to demonstrate the results of our \rl algorithm. We assume a state space $\cX\in\{0,1\}$ and action spaces $\cA^l\in\{0,1\}, \cA^f\in\{0,1\}$ for leader and the follower respectively. The state transition matrix i.e. $\tau\left(x_{t+1}|x_t,a_t\right) = .1$ if $x_t=x_{t+1}$ and $0.9$ otherwise. We assume the discounting factor $\delta = 0.6$. The rewards corresponding to the actions of one player is dependent both the state and the actions of the other player. This is tabulated in~\ref{tab:reward_state0} and~\ref{tab:reward_state1} for state $0$ and state $1$ respectively. 

\begin{figure}[!htb] 
	\centering
	\includegraphics[width=.3\textwidth]{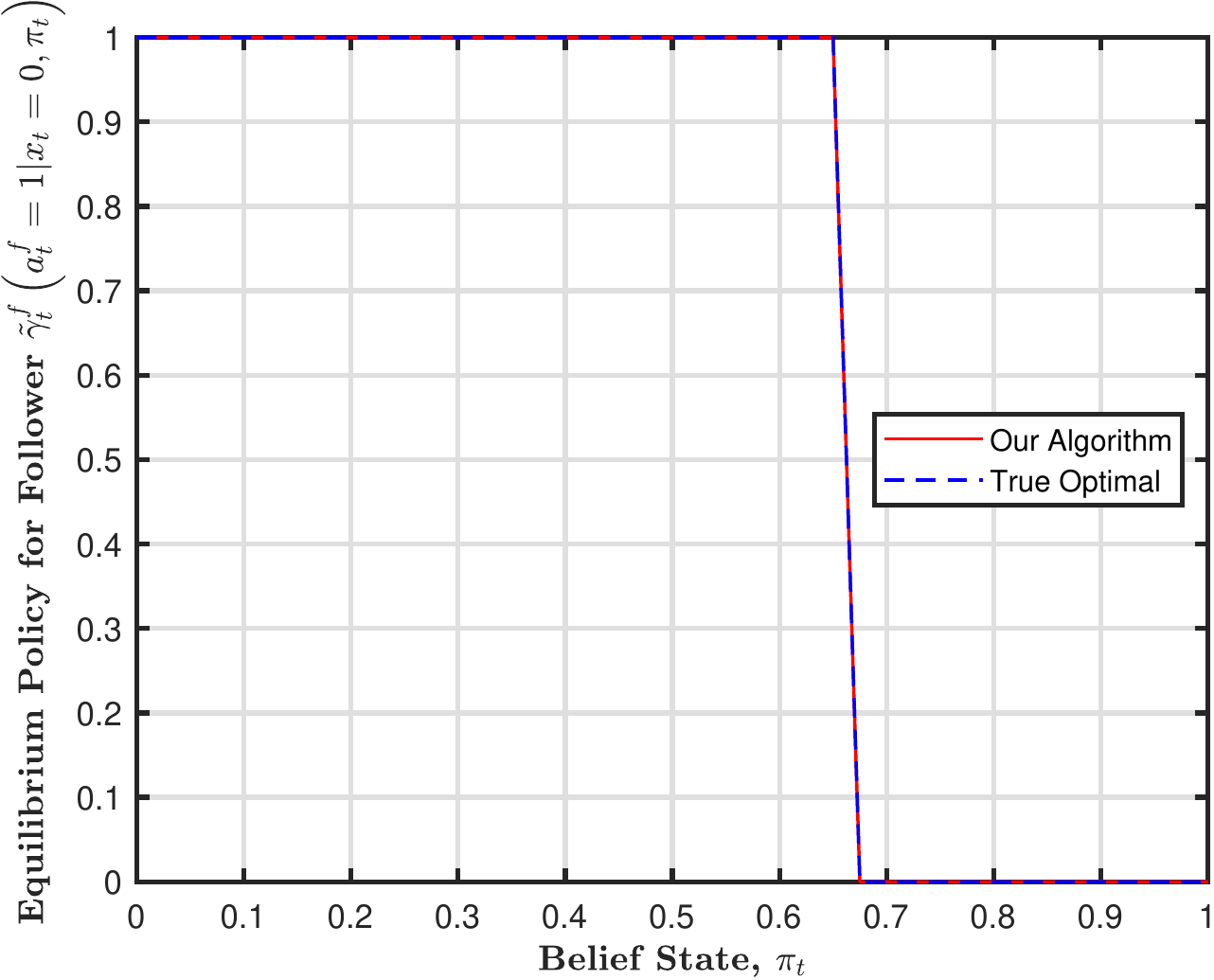} 
	\caption{$\tgamma_t^l(1|0)$: Follower's action given $x_t = 0$ for all belief state under \mpse}
	\label{fig:Figure_P1_01}
\end{figure}

\begin{figure}[!htb] 
	\centering
	\includegraphics[width=.3\textwidth]{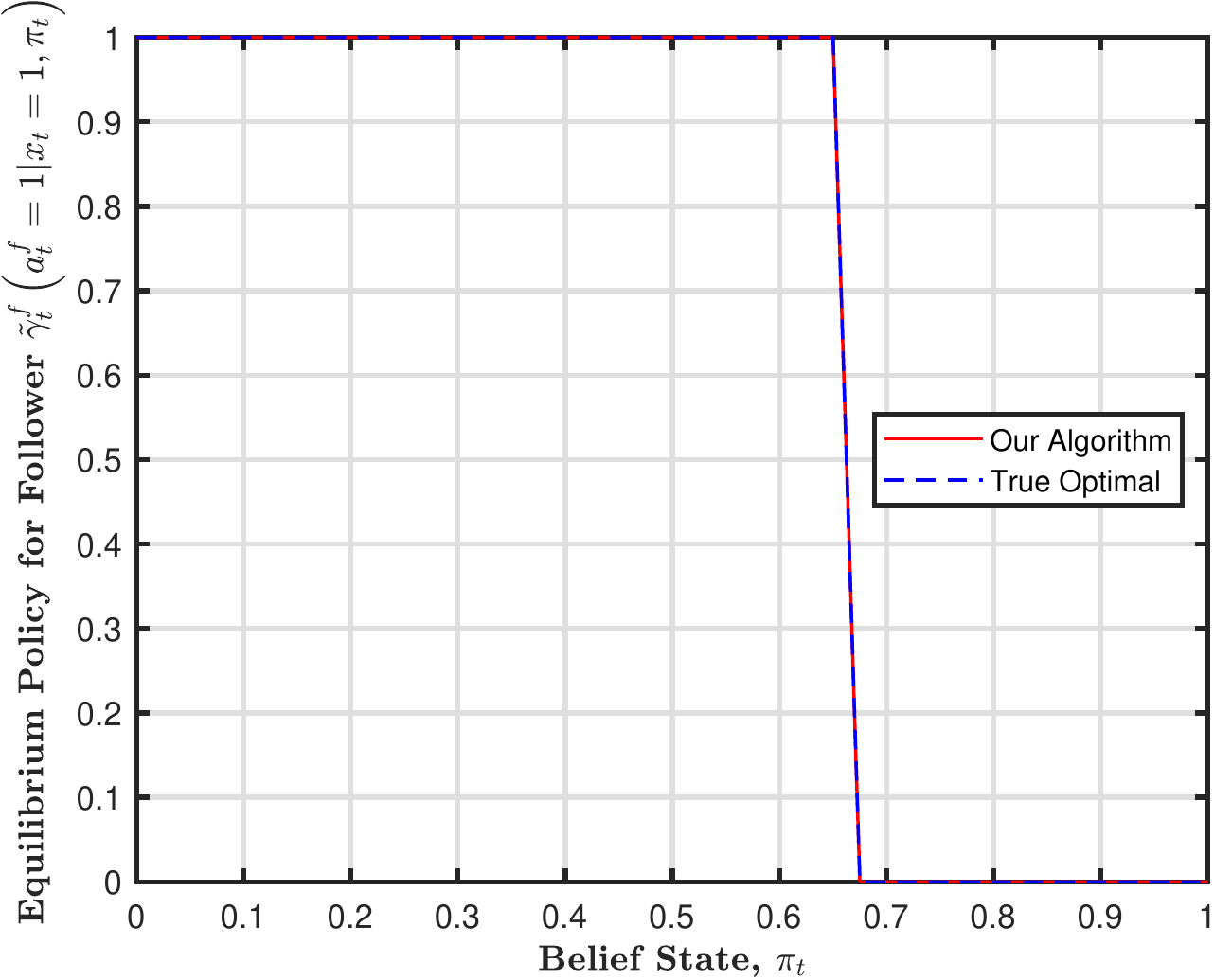} 
	\caption{$\tgamma_t^f(1|1x)$: Follower's action given $x_t = 1$ for all belief state under \mpse}
	\label{fig:Figure_P1_11}
\end{figure}

\begin{figure}[!htb] 
	\centering
	\includegraphics[width=.3\textwidth]{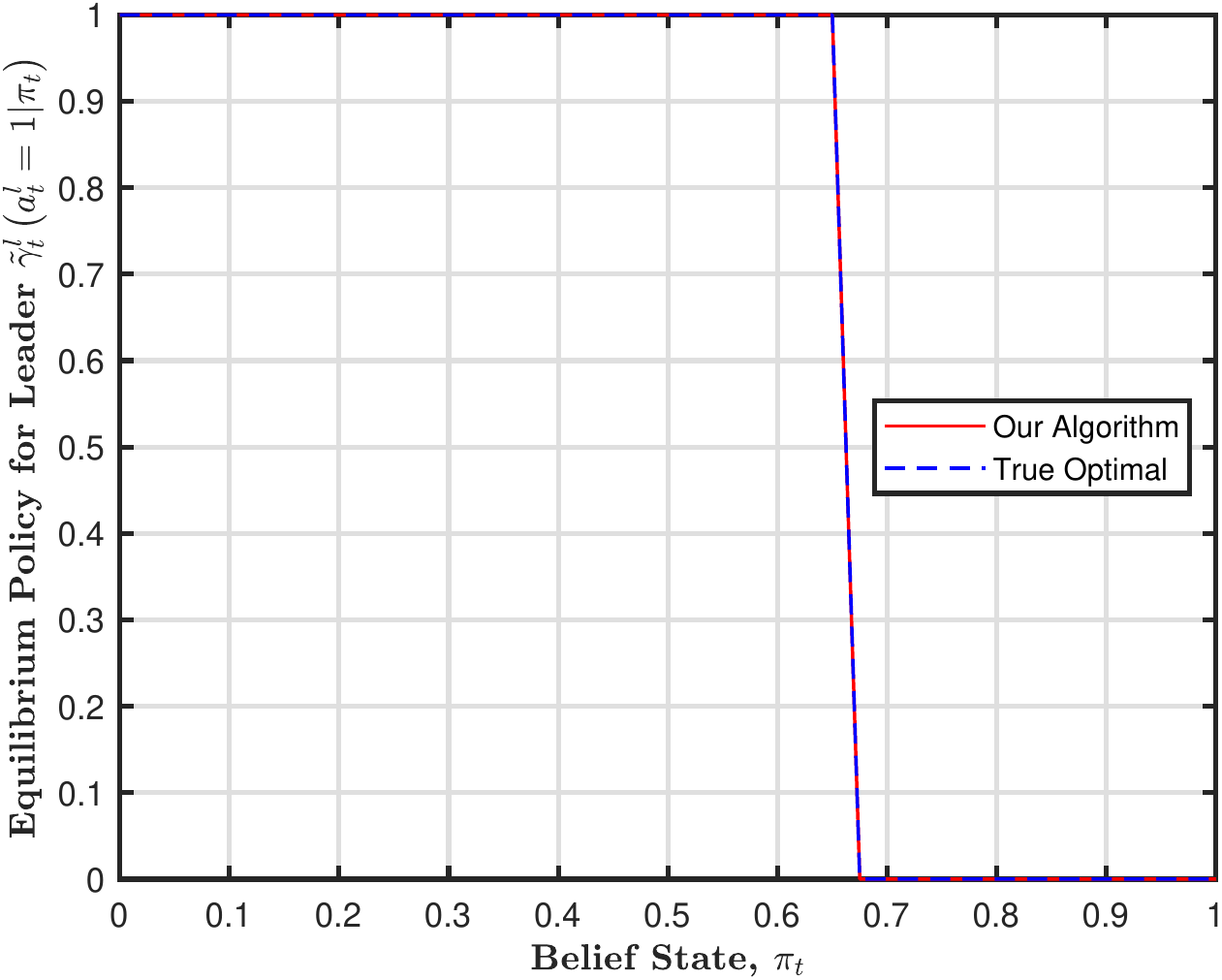} 
	\caption{$\tgamma_t^l(1|0)$: Leader's action corresponding to its \mpse policy}
	\label{fig:Figure_P2_1}
\end{figure}

\figref{fig:Figure_P1_01} and \figref{fig:Figure_P1_11} show the \mpse policies $\tgammatf$ of the follower at different values of the estimated belief state $\pih_t$ for states $x_t=0$ and $x_t=1$ respectively. The plotted graphs are the probabilities with which we choose action $a_t^i =1$.  \figref{fig:Figure_P2_1} show the \mpse policies of the leader pertaining to the probabilities of taking action $a_t^l=1$ given a belief state $\pih_t$. The plots of our algorithm are compared across the true strategy that was obtained by assuming the knowledge of the dynamics of \mdp and then solving the fixed point equation. The strategies estimated using the proposed \rl algorithm coincides with the true optimal establishing the accuracy of our algorithm.
\begin{table}[ht]
	\centering
	\caption{Reward Matrix for $x_t = 0$}
	\begin{tabular}{|cc|c|c|}
		\hline
		&      & Attacker & Attacker\\
		&      & $A1$     & $A2$ \\
		\hline
		Defender  	& $D1$ & $(2,1)$  & $(4,0)$ \\
		\hline          
		Defender    & $D2$ & $(1,0)$  & $(3,2)$ \\
		\hline
	\end{tabular}
	\label{tab:reward_state0}
\end{table}

\begin{table}[ht]
	\centering
	\caption{Reward Matrix for $x_t = 1$}
	\begin{tabular}{|cc|c|c|}
		\hline
		&      & Attacker & Attacker\\
		&      & $A1$     & $A2$ \\
		\hline
		Defender  	& $D1$ & $(2,1)$  & $(4,0)$ \\
		\hline          
		Defender    & $D2$ & $(1,0)$  & $(3,2)$ \\
		\hline
	\end{tabular}
	\label{tab:reward_state1}
\end{table}
	\section{Conclusion}
	\label{sec:Conclusion}
		\label{sec:concl}
In this paper, we analyzed a stochastic Stackelberg game where there is a leader and a follower. Follower has a private type which evolves in a controlled Markovian fashion, whose statistics are not known to both the players. We proposed a \rl algorithm, based on Expected Sarsa along with particle filters to learn the dynamics of the model by sampling and estimate the $Q_t^f$ and $Q_t^l$ functions that captures the rewards achieved following different policies from a particular state for both the players. Then, using the \mpse algorithm presented in~\cite{Vasal} to converge upon the perfect Stackelberg equilibrium of the game for both the players within a $\epsilon$ margin. 

	\medskip
	\small
	\bibliographystyle{IEEEtran}
	\bibliography{References, References1}
\end{document}